\documentclass[conference]{IEEEtran}
\IEEEoverridecommandlockouts
\usepackage{booktabs} 
\usepackage{amsmath,amssymb,amsfonts}
\usepackage[ruled,vlined]{algorithm2e}

\usepackage{algpseudocode}
\usepackage{graphicx}
\usepackage{multirow}
\usepackage{amssymb}
\usepackage{subcaption} 
\usepackage{xcolor}
\usepackage{diagbox}
\newcommand{\R}{\mathbb{R}}
\makeatletter
\renewcommand\footnoterule{%
  \kern-3\p@
  \hrule\@width.69\columnwidth
  \kern2.6\p@}
  \makeatother
\usepackage[bookmarks=false]{hyperref}
\renewcommand{\baselinestretch}{0.995}

\def\BibTeX{{\rm B\kern-.05em{\sc i\kern-.025em b}\kern-.08em
    T\kern-.1667em\lower.7ex\hbox{E}\kern-.125emX}}
\begin{document}


\title{
Exploiting Extended Krylov Subspace for the Reduction of Regular and Singular Circuit Models
}

\author{\IEEEauthorblockN{Chrysostomos~Chatzigeorgiou\textsuperscript{\textsection},~Dimitrios~Garyfallou\textsuperscript{\textsection},~George~Floros\textsuperscript{\textsection},\\~Nestor~Evmorfopoulos,~and~George~Stamoulis}

\IEEEauthorblockA{Department of Electrical and Computer Engineering, 
University of Thessaly, Volos, Greece \\ 
\{hrhatzig, digaryfa, gefloros, nestevmo, georges\}@e-ce.uth.gr }
}

\maketitle
\begingroup\renewcommand\thefootnote{\textsection}
\footnotetext{These authors contributed equally to this work}
\endgroup
\begin{abstract}
During the past decade, Model Order Reduction (MOR) has become key enabler for the efficient simulation of large circuit models. MOR techniques based on moment-matching are well established due to their simplicity and computational performance in the reduction process. 
	However, moment-matching methods based on the ordinary Krylov subspace are usually inadequate to accurately approximate the original circuit behavior.
	In this paper, we present a moment-matching method which is based on the extended Krylov subspace and exploits the superposition property in order to deal with many terminals. The proposed method can handle large-scale regular and singular circuits and generate accurate and efficient reduced-order models for circuit simulation. Experimental results on industrial IBM power grids demonstrate that our method achieves an error reduction up to 83.69\% over a standard Krylov subspace method.
	
\end{abstract}


\section{Introduction}
The ongoing miniaturization
of modern IC devices has led to extremely complex circuits.
This results in the increase of the problems associated with the analysis and simulation of their physical models.
In particular, the performance and reliable operation of ICs are largely determined by several critical subsystems such as the power distribution network, multi-conductor interconnections, and the semiconductor substrate. The electrical models of the above subsystems are very large, consisting of hundreds of millions or billions of electrical elements (mostly resistors \textit{R}, capacitors \textit{C}, and inductors \textit{L}), and their simulation 
is becoming a challenging numerical problem. Although their individual simulation is feasible, it is completely impossible to combine them 
and simulate the entire IC in many time-steps or frequencies. However, for the above subsystems it is often not necessary to fully simulate all internal state variables (node voltages and branch currents), as we only need to calculate the responses in the time or frequency domain for a small subset of output terminals (ports) and given excitations at some input ports. In these cases, the very large electrical model can be replaced by a much smaller model whose behavior at the input/output ports is similar to the behavior of the original model. This process is called Model Order Reduction (MOR).

MOR methods are divided into two main categories.
System theoretic techniques, such as Balanced Truncation (BT) \cite{gtbr}, provide very satisfactory and reliable bounds for the approximation error. However, BT techniques require the solution of Lyapunov matrix equations which are very computationally expensive, and also involve storage of dense matrices, even if the system matrices are sparse.
On the other hand, moment-matching (MM) techniques \cite{prima} are well established due to their computational efficiency in producing reduced-order models. Their drawback is that the reduced-order model depends only on the quality of the Krylov subspace.

The majority of MM methods exploit the standard or the rational Krylov subspace in order to approximate the original model. 
Authors in~\cite{morpg1,morpg2} employ rational Krylov MM methods to reduce power delivery networks. Using this projection subspace requires a heuristic and expensive parameter selection procedure, while the approximation quality is usually very sensitive to an inaccurate selection of these parameters. Moreover, in~\cite{prima,date-sing} a standard Krylov subspace is employed for the reduction of regular and singular systems, respectively. Generally, established MM methods construct the subspace only for positive directions, usually leading to a large approximated subspace to obtain a satisfactory error. Recent developments in a wide range of applications have shown that the  approximation quality of the Extended Krylov Subspace (EKS) outperforms the one of the standard Krylov subspace \cite{kzimer}. However, the application of EKS in the context of circuit simulation is not trivial.
In several problems, EKS computation involves singular circuit models and 
dense matrix manipulations, which can hinder the applicability of this subspace. 

In this paper, we introduce an EKS Moment-Matching (EKS-MM) method that greatly decreases the error induced by MM methods by approximating both ends of the spectrum. To enable the simulation of many-port models, the proposed method exploits the superposition property.
More specifically, we develop a procedure for applying EKS-MM to large-scale regular and singular models, by implementing computationally efficient transformations in order to preserve the original form of the sparse input matrices. Finally, we evaluate our methodology on industrial IBM power grids.
The rest of the paper is organized as follows. Section~\ref{MM} presents the theoretical background of MM methods for the reduction of regular and singular circuit models. Section~\ref{EKS} presents our main contributions on the application of EKS to MM methods, as well as its efficient implementation by sparse matrix manipulations for both regular and singular circuit models. Section~\ref{EXP} presents the experimental results, while conclusions are drawn in Section~\ref{CONCL}.

\section{Background} \label{MM}
\subsection{MOR by Moment-Matching}

Consider the Modified Nodal Analysis (MNA) description of an n-node, m-branch (inductive), p-input, and q-output RLC circuit in the time domain:

\begin{equation} \label{mna}
\begin{aligned}
\begin{pmatrix}
\mathbf{G} & \mathbf{W} \\
\mathbf{-W}^T & \mathbf{0} 
\end{pmatrix}\begin{pmatrix}
\mathbf{v}(t) \\
\mathbf{i}(t) 
\end{pmatrix}+\begin{pmatrix}
\mathbf{C} & \mathbf{0} \\
\mathbf{0} & \mathbf{M} 
\end{pmatrix}\begin{pmatrix}
\dot{\mathbf{v}}(t) \\
\dot{\mathbf{i}}(t) 
\end{pmatrix} = \begin{pmatrix}
\mathbf{B}_1 \\
\mathbf{0} 
\end{pmatrix}\mathbf{u}(t)\\
\mathbf{y}(t) = \begin{pmatrix}
\mathbf{L}_1 \quad
\mathbf{0} 
\end{pmatrix}\begin{pmatrix}
\mathbf{v}(t) \\
\mathbf{i}(t) 
\end{pmatrix} + \mathbf{Du}(t)
\end{aligned}
\end{equation}
where $\mathbf{G} \in\R^{n\times n}$ (node conductance matrix),  $\mathbf{C}\in\R^{n\times n}$ (node capacitance matrix), $\mathbf{M}\in\R^{m\times m}$ (branch inductance matrix), $\mathbf{W}\in\R^{n \times m}$ (node-to-branch incidence matrix), $\mathbf{v}\in\R^{n}$ (vector of node voltages), $\mathbf{i}\in\R^{m}$ (vector of inductive branch currents), $\mathbf{u} \in \R^{p} $ (vector of input excitations from current sources), $\mathbf{B}_1\in\R^{n\times p}$ (input-to-node connectivity matrix), $\mathbf{y}\in\R^{q}$ (vector of output measurements), $\mathbf{L}_1\in\R^{q\times n}$ (node-to-output connectivity matrix), $\mathbf{D} \in \R^{q\times p}$ (input-to-output connectivity matrix). Without loss of generality, in the above we assume that any voltage sources have been transformed to Norton-equivalent current sources, and that all outputs are obtained at the nodes as node voltages. Furthermore, $\dot{\mathbf{v}}(t) \equiv \frac{d\mathbf{v}(t)}{dt}$ and $\dot{\mathbf{i}}(t) \equiv \frac{d\mathbf{i}(t)}{dt}$.\\ 
If we now denote the model \textit{order} as $N \equiv n + m$, the \textit{state} vector as $\mathbf{x}(t) \equiv \begin{pmatrix}
\mathbf{v}(t) \\
\mathbf{i}(t) 
\end{pmatrix}$, and also:
\begin{equation*}
\begin{aligned}
\mathbf{A}\equiv -\begin{pmatrix}
\mathbf{G} & \mathbf{W} \\
\mathbf{-W}^T & \mathbf{0} 
\end{pmatrix},\quad  \mathbf{E} \equiv \begin{pmatrix}
\mathbf{C} & \mathbf{0} \\
\mathbf{0} & \mathbf{M} 
\end{pmatrix},\\ \mathbf{B}\equiv \begin{pmatrix}
\mathbf{B}_1 \\
\mathbf{0} 
\end{pmatrix}, \quad \mathbf{L}\equiv \begin{pmatrix}
\mathbf{L}_1 \quad
\mathbf{0} 
\end{pmatrix}
\end{aligned}
\end{equation*}
then expression (\ref{mna}) can be written in the following generalized state-space form, or so-called \textit{descriptor} form:

\begin{equation}
\begin{aligned} \label{state}
\mathbf{E}\frac{d \mathbf{x}(t)}{d t} = \mathbf{A x}(t) + \mathbf{Bu}(t), \quad
\mathbf{y}(t) = \mathbf{L x}(t) + \mathbf{Du}(t)
\end{aligned}
\end{equation}
The objective of MOR is to produce a reduced-order model:

\begin{equation}
\begin{aligned}\label{state_red}
\mathbf{\tilde E} \frac{d \mathbf{\tilde x}(t)}{d t} =\mathbf{\tilde A} \mathbf{\tilde x}(t) + \mathbf{\tilde B} \mathbf{u(t)}, \quad
\mathbf{\tilde y}(t) = \mathbf{\tilde L \tilde x}(t) + \mathbf{Du}(t)
\end{aligned}
\end{equation}				
where  $\mathbf{\tilde A}, \mathbf{\tilde E} \in \R^{r\times r} $, $\mathbf{\tilde B} \in \R^{r\times p} $, $\mathbf{\tilde L} \in \R^{q\times r} $. The reduced model has order $r<<N$, and the output error is bounded as $||\mathbf{\tilde{y} }(t) -\mathbf{y}(t)||_2 < \varepsilon||\mathbf{u}(t)||_2$ for given input $\mathbf{u}(t)$ and given small $\varepsilon$. The bound in the output error can be equivalently written in the frequency domain as $||\mathbf{\tilde{y} }(s) -\mathbf{y}(s)||_2 < \varepsilon||\mathbf{u}(s)||_2$ via the Plancherel's theorem \cite{plancherel}. If

\begin{equation*}
\begin{aligned}
\mathbf{H}(s) = \mathbf{L}(s\mathbf{E} - \mathbf{A})^{-1} \mathbf{B} + \mathbf{D}\\
\mathbf{\tilde H}(s) =  \mathbf{\tilde L}(s\mathbf{\tilde E} - \mathbf{\tilde A})^{-1} \mathbf{\tilde B} + \mathbf{D}
\end{aligned}
\end{equation*}
are the transfer functions of the original and the reduced-order model, respectively, then the output error in frequency domain is:

\begin{equation}
\begin{aligned}
||\mathbf{\tilde{y} }(s) -\mathbf{y}(s)||_2 = ||\mathbf{\tilde{H}}(s) \mathbf{u}(s) - \mathbf{H}(s)\mathbf{u}(s)||_2 \\ \leq ||\mathbf{\tilde{H}}(s) - \mathbf{H}(s)||_\infty||\mathbf{u}(s)||_2	
\end{aligned}
\end{equation}
where $||.||_\infty$ is the induced $\mathcal{L}_2$ matrix norm, or $\mathcal{H}_\infty$ norm, of a rational transfer function. Therefore, in order to bound the output error, we need to bound the distance between the transfer functions $||\mathbf{\tilde{H}}(s) - \mathbf{H}(s)||_\infty < \varepsilon$.

The most important and successful MOR methods for linear systems are based on MM. They are very efficient in circuit simulation problems and are formulated 
in a way that has a direct application to the linear model of (\ref{state}). 

By applying the Laplace transform to (\ref{state}), we obtain the $s$ domain equations as:

\begin{equation}
\begin{aligned}
s\mathbf{EX}(s)- \mathbf{X}(0) = \mathbf{A}\mathbf{X}(s) + \mathbf{B}\mathbf{U}(s) \\
\mathbf{Y}(s) = \mathbf{L}\mathbf{X}(s) + \mathbf{D}\mathbf{U}(s)
\end{aligned}
\end{equation}
Assuming that $\mathbf{X}(0) = 0$ and that an impulse response is applied to $\mathbf{U}(s)$ (i.e. $\mathbf{U}(s) = 1$), then the above system of equations can be written as follows:
\begin{equation}
\begin{aligned}
(s\mathbf{E}- \mathbf{A} )\mathbf{X}(s)=\mathbf{B}, \quad
\mathbf{Y}(s) = \mathbf{L}\mathbf{X}(s) + \mathbf{D}
\end{aligned}
\end{equation}
and by expanding the Taylor series of $\mathbf{X}(s)$ around zero, we derive the below equation:
\begin{equation}
(s\mathbf{E}- \mathbf{A})(\mathbf{x}_0 + \mathbf{x}_1s + \mathbf{x}_2s^2 + \dots) = \mathbf{B}
\end{equation} 
The transfer function of (\ref{state}) is a function of $s$, and can be expanded into a moment expansion around $s = 0$ as follows:
\begin{equation} \label{mom2}
\mathbf{H}(s) = \mathbf{M}_0 + \mathbf{M}_1s + \mathbf{M}_2s^2 + \mathbf{M}_3s^3 \dots 
\end{equation}
where $\mathbf{M}_0$, $\mathbf{M}_1$, $\mathbf{M}_2$, $\mathbf{M}_3$, $\dots$ are the moments of the transfer function. Specifically, in circuit simulation problems, $\mathbf{M}_0$ is the DC solution of the linear system. This means that the inductors of the circuit are considered as short
circuits, and the capacitors as open circuits. Moreover, $\mathbf{M}_1$  is the Elmore delay of the linear model, which is defined as the time required for a signal at the input port to reach the output port. 
Finally, $\mathbf{M}_i$ is related to the system matrices as:

\begin{equation}
\mathbf{M}_i = \mathbf{L} (\mathbf{A}^{-1}\mathbf{E})^i\mathbf{A}^{-1}\mathbf{B}
\end{equation}
The goal of MM reduction techniques  is the derivation of a reduced-order model where some moments $\mathbf{\tilde M}_i$ of the reduced-order transfer function $\mathbf{\tilde{H}}(s)$ match some moments of the original transfer function $\mathbf{H}(s)$.

Let us now denote the two projection matrices onto a lower dimensional subspace as $\mathbf{W} \in \R^{N \times r} $ and $\mathbf{V}\in \R^{r \times N}$, respectively. These matrices can be derived from the associated moment vectors using one or more expansion points. As a result, if we assume that $s=0$, then the matrices $\mathbf{W}$ and $\mathbf{V}$ are defined as follows: 
\begin{equation}
\begin{aligned}
range(\mathbf{W})= span\{\mathbf{B},(\mathbf{A}^{-1}\mathbf{E})\mathbf{B},\dots, (\mathbf{A}^{-1}\mathbf{E})^r\mathbf{B}\}\\
range(\mathbf{V}) = span\{\mathbf{L},(\mathbf{A}^{-1}\mathbf{E})^{-T}\mathbf{L},\dots, (\mathbf{A}^{-T}\mathbf{E}^T)^r\mathbf{L}\}
\end{aligned}
\end{equation}
The computed reduced-order model matches the first $2r$ moments and is obtained by the following matrices:
\begin{equation}
\mathbf{\tilde E} = \mathbf{W}^T\mathbf{ E}\mathbf{V}, \quad \mathbf{\tilde A} = \mathbf{W}^T\mathbf{ A}\mathbf{V}, \quad \mathbf{\tilde B} = \mathbf{W}^T\mathbf{ B}, \quad \mathbf{\tilde L} = \mathbf{ L}\mathbf{V}
\end{equation}
This reduced model provides a good approximation around the DC point. Finally, in case we employ an one-sided Krylov method, which is usually the case, the matrix $\mathbf{W}$ can be set equal to $\mathbf{V}$, an equality that also holds for symmetric systems.

\subsection{Handling of Singular Descriptor Models} \label{Singular}

In certain circuit simulation problems, the matrix $\mathbf{E}$ might be singular. A method for dealing with such models is to compute spectral projections onto the left and right deflating subspaces corresponding to the finite eigenvalues of the model, which is computationally prohibitive for large-scale systems. However, singular descriptor models typically result when there are some nodes, say $n_2$,  where no capacitance is connected, leading to corresponding all-zero rows and columns in the submatrix $\mathbf{C}$. Note that in case the circuit contains no voltage sources, the submatrix $\mathbf{M}$ of inductive branches is always nonsingular. If the $n_2$ nodes with no capacitance connection are enumerated last, and the remaining $n_1 = n -n_2$ nodes first, then (\ref{mna}) can be partitioned as follows:

\begin{equation} \label{12}
\begin{aligned}
\begin{pmatrix}
\mathbf{G}_{11} &\mathbf{G}_{12} & \mathbf{W}_1 \\
\mathbf{G}_{12}^T &\mathbf{G}_{22} & \mathbf{W}_2 \\
-\mathbf{W}_1^T & -\mathbf{W}_2^T&\mathbf{0} 
\end{pmatrix}\begin{pmatrix}
\mathbf{v}_1(t) \\
\mathbf{v}_2(t) \\
\mathbf{i}(t) 
\end{pmatrix}+\\\begin{pmatrix}
\mathbf{C}_{1} & \mathbf{0} & \mathbf{0} \\
\mathbf{0} & \mathbf{0} & \mathbf{0} \\
\mathbf{0} &\mathbf{0}& \mathbf{M} 
\end{pmatrix}\begin{pmatrix}
\dot{\mathbf{v}}_1(t) \\
\dot{\mathbf{v}}_2(t) \\
\dot{\mathbf{i}}(t) 
\end{pmatrix} = \begin{pmatrix}
\mathbf{B}_1 \\
\mathbf{B}_2 \\
\mathbf{0} 
\end{pmatrix}\mathbf{u}(t)\\
\mathbf{y}(t) = \begin{pmatrix}
\mathbf{L}_1\quad
\mathbf{L}_2 \quad 
\mathbf{0} 
\end{pmatrix}\begin{pmatrix}
\mathbf{v}_1(t) \\
\mathbf{v}_2(t) \\
\mathbf{i}(t) 
\end{pmatrix} + \mathbf{Du}(t)
\end{aligned}
\end{equation}
where $\mathbf{G}_{11} \in \R^{n_1\times n_1}$, $\mathbf{G}_{12} \in \R^{n_1\times n_2}$, $\mathbf{G}_{22} \in \R^{n_2\times n_2}$, $\mathbf{W}_{1} \in \R^{n_1\times m}$, $\mathbf{W}_{2} \in \R^{n_2\times m}$, $\mathbf{C}_{1} \in \R^{n_1\times n_1}$, $\mathbf{v}_{1} \in \R^{n_1}$, $\mathbf{v}_{2} \in \R^{n_1}$, $\mathbf{B}_{1} \in \R^{n_1 \times p}$, $\mathbf{B}_{2} \in \R^{n_2 \times p}$, $\mathbf{L}_{1} \in \R^{q \times n_1}$, and $\mathbf{L}_{2} \in \R^{q\times n_2}$.\\
Assuming now that the submatrix $\mathbf{G}_{22}$ is nonsingular 
(a sufficient condition for this is at least one resistive connection from any of the $n_2$ non-capacitive nodes to ground), the second row of (\ref{12}) can be solved for $\mathbf{v}_2(t)$ as follows: 
\begin{equation}\label{13}
\begin{aligned}
\mathbf{v}_2(t) =\mathbf{G}_{22}^{-1}\mathbf{B}_2\mathbf{u}(t) - \mathbf{G}_{22}^{-1}\mathbf{G}_{12}^T\mathbf{v}_1(t) - \mathbf{G}_{22}^{-1}\mathbf{W}_{2}\mathbf{i}(t) 
\end{aligned}
\end{equation}
The above can be substituted to the first and third row of (\ref{12}), as well as the output part of (\ref{12}), to give:
\begin{equation*}
\begin{aligned}
(\mathbf{G}_{11} - \mathbf{G}_{12}\mathbf{G}_{22}^{-1}\mathbf{G}_{12}^T)\mathbf{v}_1(t)+ (\mathbf{W}_1- \mathbf{G}_{12}\mathbf{G}_{22}^{-1}\mathbf{W}_{2})\mathbf{i}(t)\\+\mathbf{C}_{1}\dot{\mathbf{v}}_1(t)= (\mathbf{B}_1-\mathbf{G}_{12}\mathbf{G}_{22}^{-1}\mathbf{B}_2)\mathbf{u}(t)
\end{aligned}
\end{equation*}
\begin{equation*}
\begin{aligned}
(\mathbf{W}_{2}^T\mathbf{G}_{22}^{-1}\mathbf{G}_{12}^T -\mathbf{W}_{1}^T)\mathbf{v}_1(t)+ \mathbf{W}_{2}^T\mathbf{G}_{22}^{-1}\mathbf{W}_{2}\mathbf{i}(t)+\mathbf{M}\dot{\mathbf{i}}(t)\\= \mathbf{W}_{2}^T\mathbf{G}_{22}^{-1}\mathbf{B}_{2}\mathbf{u}(t)
\end{aligned}
\end{equation*}
\begin{equation*}
\begin{aligned}
\mathbf{y}(t) = (\mathbf{L}_{1}- \mathbf{L}_{2}\mathbf{G}_{22}^{-1}\mathbf{G}_{12}^T)\mathbf{v}_{1}(t) - \mathbf{L}_{2}\mathbf{G}_{22}^{-1}\mathbf{W}_{2}\mathbf{i}(t)\\ + (\mathbf{L}_{2}\mathbf{G}_{22}^{-1}\mathbf{B}_{2}+\mathbf{D})\mathbf{u}(t)
\end{aligned}
\end{equation*} 
This can be put together in the following descriptor form:
\begin{equation} \label{111}
\begin{aligned}
\begin{pmatrix}
\mathbf{C}_{1} & \mathbf{0} \\
\mathbf{0} & \mathbf{M} 
\end{pmatrix}\begin{pmatrix}
\dot{\mathbf{v}}_1(t) \\
\dot{\mathbf{i}}(t) 
\end{pmatrix} = \\-\begin{pmatrix}
\mathbf{G}_{11} - \mathbf{G}_{12}\mathbf{G}_{22}^{-1}\mathbf{G}_{12}^T & \mathbf{W}_1- \mathbf{G}_{12}\mathbf{G}_{22}^{-1}\mathbf{W}_{2} \\
\mathbf{W}_{2}^T\mathbf{G}_{22}^{-1}\mathbf{G}_{12}^T -\mathbf{W}_{1}^T & \mathbf{W}_{2}^T\mathbf{G}_{22}^{-1}\mathbf{W}_{2}
\end{pmatrix}\begin{pmatrix}
\mathbf{v}_1(t) \\
\mathbf{i}(t) 
\end{pmatrix} \\+ \begin{pmatrix}
\mathbf{\mathbf{B}_1-\mathbf{G}_{12}\mathbf{G}_{22}^{-1}\mathbf{B}_2} \\
\mathbf{W}_{2}^T\mathbf{G}_{22}^{-1}\mathbf{B}_{2} 
\end{pmatrix}\mathbf{u}(t)\\
\mathbf{y}(t)= \begin{pmatrix}
\mathbf{L}_{1}- \mathbf{L}_{2}\mathbf{G}_{22}^{-1}\mathbf{G}_{12}^T \quad
\mathbf{L}_{2}\mathbf{G}_{22}^{-1}\mathbf{W}_{2} 
\end{pmatrix}\begin{pmatrix}
{\mathbf{v}_1}(t) \\
{\mathbf{i}}(t) 
\end{pmatrix}\\ + (\mathbf{L}_{2}\mathbf{G}_{22}^{-1}\mathbf{B}_{2}+\mathbf{D})\mathbf{u}(t)
\end{aligned}
\end{equation}
The above is a nonsingular (i.e. regular) state-space model which can be reduced normally.

\section{Extended Krylov Subspace for MOR}\label{EKS}

\subsection{EKS Moment-Matching (EKS-MM)} \label{3.4}
The essence of MM methods is to iteratively compute a projection subspace, and then project the original system into this subspace in order to obtain the reduced-order model of (\ref{state_red}). The dimension of the projection subspace is increased
in every iteration, until an a-priory selection of the moments is matched.
More specifically, if $r$ is the desired order for the reduced system  and $k= \frac{r}{p}$ is the number of moments, then $\mathbf{X} \in \R^{N \times r}$ ($r<<N$) is a projection matrix whose columns span the $k$-dimensional Krylov subspace:
 \begin{equation*}
  \mathcal{K}_k(\mathbf{A}_{E},\mathbf{B}_{E})=  
  span  \{\mathbf{B}_{E},\mathbf{A}_{E}\mathbf{B}_{E},  \mathbf{A}_{E}^{2}\mathbf{B}_{E},\dots,\mathbf{A}_{E}^{k-1}\mathbf{B}_{E}\}
 \end{equation*}
where 
\begin{equation*}
\mathbf{A}_{E} \equiv \mathbf{A}^{-1}\mathbf{E},\quad \mathbf{B}_{E} \equiv \mathbf{A}^{-1}\mathbf{B}
\end{equation*}
Then, the reduced-order model is obtained through the following matrix transformations:
\begin{equation} \label{20}
\mathbf{\tilde E} = \mathbf{X}^T\mathbf{ E}\mathbf{X}, \quad \mathbf{\tilde A} = \mathbf{X}^T\mathbf{ A}\mathbf{X}, \quad \mathbf{\tilde B} = \mathbf{X}^T\mathbf{ B}, \quad \mathbf{\tilde L} = \mathbf{ L}\mathbf{X}
\end{equation}
with $\mathbf{\tilde A}, \mathbf{\tilde E} \in \R^{r\times r} $, $\mathbf{\tilde B} \in \R^{r\times p} $, $\mathbf{\tilde L} \in \R^{q\times r} $.

The projection process is independent of the subspace selection, but its effectiveness is critically dependent on the chosen
subspace. As a result, one choice is to consider the rational Krylov subspace~\cite{morpg1,morpg2}.  However, this projection subspace requires the input of a number of shift parameters, whose choice greatly affects the produced reduced-order model. 
The reason for this is that it relies on unclear heuristics and is highly problem-dependent. In order to address this issue, the standard Krylov subspace ~\cite{prima,date-sing} $\mathcal{K}_k(\mathbf{A}_{E},\mathbf{B}_{E}) $ must be enriched with information from the subspace $\mathcal{K}_k(\mathbf{A}_{E}^{-1},\mathbf{B}_{E})$, which corresponds to the inverse matrix $\mathbf{A}_{E}^{-1}$, leading to EKS:
\begin{equation*}
\label{Eq:eksm}
\mathcal{K}_k^E(\mathbf{A}_{E},\mathbf{B}_{E}) = \mathcal{K}_k(\mathbf{A}_{E},\mathbf{B}_{E}) + \mathcal{K}_k(\mathbf{A}_{E}^{-1},\mathbf{B}_{E}) =
\end{equation*}
\begin{equation}\label{19}
span \{\mathbf{B}_{E},  \mathbf{A}_{E}^{-1}\mathbf{B}_{E}, \mathbf{A}_{E}\mathbf{B}_{E},\mathbf{A}_{E}^{-2}\mathbf{B}_{E}, \mathbf{A}_{E}^{2}\mathbf{B}_{E},\dots,  
\end{equation}
\begin{equation*}
\mathbf{A}_{E}^{-(k-1)}\mathbf{B}_{E}, \mathbf{A}_{E}^{k-1}\mathbf{B}_{E}\}
\end{equation*} 
The Arnoldi procedure~\cite{mat} that computes EKS begins with the pair $\{ \mathbf{B}_{E}, \mathbf{A}_{E}^{-1}\mathbf{B}_{E} \}$, and then generates a sequence of extended subspaces $\mathcal{K}_k^E(\mathbf{A}_{E},\mathbf{B}_{E})$ in order to compute the matrix $\mathbf{X}\in \R^{N \times 2r}$ and produce the reduced-order model as described in (\ref{20}). EKS can be considered a special case of the rational Krylov subspace with two expansion points, one expansion point at zero and one at infinity.
The complete procedure is given in Algorithm \ref{eksm_alg}.


\LinesNumbered

\begin{algorithm}[hbt!]
\footnotesize
\SetAlgoLined
\DontPrintSemicolon
\KwIn{$ \mathbf{A}_{E} \equiv \mathbf{A^{-1}E}, \mathbf{B}_{E} \equiv \mathbf{A^{-1}B}$, desired order $r$, \#ports $p$}    
\KwOut{$ \mathbf{X} $}

\SetKwFunction{FMain}{\texttt{compute\_EKS}}
\SetKwProg{Fn}{Function}{:}{}
\Fn{\FMain{$\mathbf{A}_{E},\mathbf{B}_{E},r$}}{

$j=1$\\
$\mathbf{X}^{(j)} = \texttt{qr}([\mathbf{B}_{E}, \mathbf{A}_{E}^{-1}\mathbf{B}_{E}])$\\
$k= \frac{r}{p}$\\

\While{ $(j<k)$ } {
$k_1=2p(j-1)$; $k_2 = k_1+p$; $k_3 = 2pj$ \\
$\mathbf{X}_1  = [\mathbf{A}_{E}\mathbf{X}^{(j)}(:,k_1+1:k_2),\mathbf{A}_{E}^{-1}\mathbf{X}^{(j)}(:,k_2+1:k_3)]$\\
$\mathbf{X}_2 =$ \texttt{orth\_wrt}$(\mathbf{X}_1, \mathbf{X}^{(j)}, p) $\\
$\mathbf{X}_3 = \texttt{qr}(\mathbf{X}_2) $ \\
$\mathbf{X}^{(j+1)} = [\mathbf{X}^{(j)},\mathbf{X}_3]$\quad\\
$j=j+1$\\
}
$\mathbf{X} = \mathbf{X}(:,1:2r)$\\
\textbf{return} $\mathbf{X}$\\
}
\textbf{End Function}

\caption{EKS computation by Arnoldi procedure}%
\label{eksm_alg}
\end{algorithm}
\normalsize
At this point, we can elaborate on some aspects regarding the efficient implementation of the proposed EKS procedure:
\subsubsection{\textbf{Sparse matrix inputs}} 
It is worth mentioning that Algorithm \ref{eksm_alg} does not require matrices $\mathbf{A}_{E} \equiv \mathbf{A}^{-1}\mathbf{E}$, $\mathbf{B}_{E}\equiv \mathbf{A}^{-1}\mathbf{B}$ as inputs, but only the sparse system matrices $\mathbf{A}$, $\mathbf{E}$ are necessary. This is due to the fact that the generally dense inverse matrices are only needed in products with $p$ vectors (initially in step 3) and $2pj$ vectors (in step 7 at every iteration, where the iteration count $j$ is typically very small and thus $2pj<<N$). 
These products can be implemented as sparse linear solves ($\mathbf{E}\mathbf{Y} = \mathbf{R}$ and $\mathbf{A}\mathbf{Y} = \mathbf{R}$) by employing any sparse direct~\cite{dir} or iterative~\cite{CMG} algorithm.
\subsubsection{\textbf{Orthogonalization in steps 3 and 9}} A modified Gram-Schmidt procedure \cite{mat} is employed to implement the corresponding $\texttt{qr()}$ procedures.
\subsubsection{\textbf{Orthogonalization in step 8}}
	In order to perform orthogonalization with respect to matrix $\mathbf{X}^{(j)}$, we employ the following Gram-Schmidt procedure \cite{mat} as shown in Algorithm \ref{orth}.
\begin{algorithm}[hbt!]
\footnotesize
\SetAlgoLined
\DontPrintSemicolon
\KwIn{$ \mathbf{X}_1,\mathbf{X}^{(j)}$, \#ports p }    
\KwOut{$ \mathbf{X}_2 $}

\SetKwFunction{FMain}{orth\_wrt}
\SetKwProg{Fn}{Function}{:}{}
\Fn{\FMain{$\mathbf{X}_1,\mathbf{X}^{(j)},p$}}{

\For{$k_1=1,\dots,j$}{
$k_2=2p(k_1-1)$; $k_3=2pk_1;$\\
$\mathbf{X}_2 = \mathbf{X}_1-\mathbf{X}^{(j)}(:,k_2+1:k_3)\mathbf{X}^{(j)T}(:,k_2+1:k_3)\mathbf{X}_1$
}
\textbf{return} $\mathbf{X}_2$\\
}
\textbf{End Function}

\caption{Orthogonalization w.r.t. another matrix}%
\label{orth}
\end{algorithm}
\normalsize
\subsection{Sparse Implementation for Singular Descriptor Models} \label{5.3}
Algorithm \ref{eksm_alg} is computationally inefficient for the reduction of the model given in (\ref{111}), which results from the regularization of a singular descriptor model, since the inversion of $\mathbf{G}_{22}$ renders the matrices dense and hinders the solution procedure. In this subsection, we present efficient ways to implement the EKS algorithm by preserving the original sparse form of the system matrices.
\subsubsection{\textbf{Construction of RHS}} The input-to-state and state-to-output connectivity matrices
\begin{equation} \label{421}
\mathbf{B} \equiv \begin{pmatrix}
\mathbf{B}_1-\mathbf{G}_{12}\mathbf{G}_{22}^{-1}\mathbf{B}_2 \\
\mathbf{W}_{2}^T\mathbf{G}_{22}^{-1}\mathbf{B}_{2} 
\end{pmatrix} , \quad \mathbf{L}^T \equiv \begin{pmatrix}
\mathbf{L}_{1}^T- \mathbf{G}_{12}\mathbf{G}_{22}^{-1}\mathbf{L}_{2}^T \\
\mathbf{W}_{2}^T\mathbf{G}_{22}^{-1}\mathbf{L}_{2}^T 
\end{pmatrix} 
\end{equation} 
are explicitly constructed to compute the input matrix $\mathbf{B}_E$ of Algorithm \ref{eksm_alg}, and to obtain the reduced order model through (\ref{20}). The products $\mathbf{G}_{22}^{-1}\mathbf{B}_2$ and $\mathbf{G}_{22}^{-1}\mathbf{L}_{2}^T$ are computed by $p$ and $q$ sparse linear solves, respectively.
\subsubsection{\textbf{Sparse linear system solutions}} The system matrix
\begin{equation}\label{422}
\mathbf{A} \equiv - \begin{pmatrix}
\mathbf{G}_{11} - \mathbf{G}_{12}\mathbf{G}_{22}^{-1}\mathbf{G}_{12}^T & \mathbf{W}_1- \mathbf{G}_{12}\mathbf{G}_{22}^{-1}\mathbf{W}_{2} \\
\mathbf{W}_{2}^T\mathbf{G}_{22}^{-1}\mathbf{G}_{12}^T -\mathbf{W}_{1}^T & \mathbf{W}_{2}^T\mathbf{G}_{22}^{-1}\mathbf{W}_{2}
\end{pmatrix}
\end{equation}
of the model given in (\ref{111}) is rendered dense due to the inversion of $\mathbf{G}_{22}$. The linear system solutions with $\mathbf{A}$ in steps 3, 7 of Algorithm \ref{eksm_alg} can be handled by partitioning the RHS of these systems conformally to $\mathbf{A}$, i.e. $\mathbf{R} =  \begin{pmatrix}
\mathbf{R}_1 \\
\mathbf{R}_2 \\
\end{pmatrix}$
with $\mathbf{R}_1 \in \R^{n_1\times p}$, $\mathbf{R}_2\in \R^{m \times p}$, 
and implementing their solution efficiently by keeping all the sub-blocks in their original sparse form as follows: 
\begin{equation} \label{122}
\begin{aligned}
\begin{pmatrix}
-\mathbf{G}_{11} &-\mathbf{W}_{1} & -\mathbf{G}_{12} \\
\mathbf{W}_1^T &\mathbf{0} & \mathbf{W}_2^T \\
-\mathbf{G}_{12}^T & -\mathbf{W}_2&-\mathbf{G}_{22} 
\end{pmatrix}\begin{pmatrix}
\mathbf{X}_1 \\
\mathbf{X}_2 \\
\mathbf{T}
\end{pmatrix} = \begin{pmatrix}
\mathbf{R}_1 \\
\mathbf{R}_2 \\
\mathbf{0}
\end{pmatrix}
\end{aligned}
\end{equation}
where $\mathbf{T} \in \R^{n_2 \times p}$ is a temporary sub-matrix.
\subsubsection{\textbf{Sparse matrix-vector products}} The matrix-vector products with $\mathbf{X}^{(j)}$ in step 7 of Algorithm \ref{eksm_alg} can be implemented efficiently by observing that:
\begin{equation}\label{424}
\begin{aligned}
\mathbf{A} =\begin{pmatrix}
-\mathbf{G}_{11}  & -\mathbf{W}_{1}\\
\mathbf{W}_1^T &\mathbf{0}
\end{pmatrix}+ \begin{pmatrix}
\mathbf{G}_{12}\mathbf{G}_{22}^{-1}\mathbf{G}_{12}^T & \mathbf{G}_{12}\mathbf{G}_{22}^{-1}\mathbf{W}_{2}\\
-\mathbf{W}_{2}^T\mathbf{G}_{22}^{-1}\mathbf{G}_{12}^T&-\mathbf{W}_{2}^T\mathbf{G}_{22}^{-1}\mathbf{W}_{2}\end{pmatrix} \\=
\begin{pmatrix}
-\mathbf{G}_{11}  & -\mathbf{W}_{1}\\
\mathbf{W}_1^T &\mathbf{0}
\end{pmatrix} + \begin{pmatrix}
-\mathbf{G}_{12} \\
\mathbf{W}_{2}^T
\end{pmatrix}\mathbf{G}_{22}^{-1}\begin{pmatrix}
-\mathbf{G}_{12}^T &
-\mathbf{W}_{2}
\end{pmatrix}
\end{aligned}
\end{equation}
Therefore, the product $\mathbf{A} \mathbf{X}^{(j)}$ with  $p$ vectors $\mathbf{X}^{(j)}$ can be carried out by a sparse solve $\mathbf{G}_{22}\mathbf{X} = \begin{pmatrix}
-\mathbf{G}_{12}^T &
-\mathbf{W}_{2}
\end{pmatrix}\mathbf{K}^{(j)}
$, followed by a sum of products $\begin{pmatrix}
-\mathbf{G}_{11}  & -\mathbf{W}_{1}\\
\mathbf{W}_1^T &\mathbf{0}
\end{pmatrix}\mathbf{K}^{(j)} + \begin{pmatrix}
-\mathbf{G}_{12} \\
\mathbf{W}_{2}^T
\end{pmatrix}\mathbf{X}$.
\vspace{1.1mm}
\subsubsection{\textbf{Construction of system matrix}} In order to construct and then reduce the dense system matrix of (\ref{422}), we need to employ sparse solves with the submatrix $\mathbf{G}_{22}$. Since usually $n_2 << n_1$, it is better to first compute the left-solves $\mathbf{G}_{12}\mathbf{G}_{22}^{-1}$ and $\mathbf{W}_{2}^T\mathbf{G}_{22}^{-1}$, followed by products with $\mathbf{G}_{12}^T$ and $\mathbf{W}_{2}$. The left-solves can be performed as $\mathbf{G}_{22}\mathbf{X} =\mathbf{G}_{12}$ and $\mathbf{G}_{22}\mathbf{X} =\mathbf{W}_{2}^T$, where $\mathbf{X}$ contains the \textit{rows} of each left-solve.  

\subsection{Superposition Property of LTI Models}
While in the previous subsections we emphasized on the efficient execution of the proposed methodology, it still can not handle many-terminal models. To this end, we consider the superposition principal of LTI models. 
Using the superposition property, the output response of the initial multi-input multi-output (MIMO) descriptor model of (\ref{state}) can be computed as the sum of the output responses of the following single-input multi-output (SIMO) subsystems  as:
\begin{equation}
\begin{aligned} \label{super_state}
\mathbf{E}\frac{d \mathbf{x}(t)}{d t} = \mathbf{A x}(t) + \mathbf{B}_i\mathbf{u}_i(t), \quad
\mathbf{y}_i(t) = \mathbf{L x}(t) + \mathbf{Du}_i(t)
\end{aligned}
\end{equation}
where $\mathbf{B}_i$ is a  matrix with only one nonzero column of the input-to-node-connectivity matrix $\mathbf{B}$, and $i=1,\dots,p$. From these relations, it can be derived that $\mathbf{y}(t)= \sum_{n=1}^{p}\mathbf{y}_i(t)$ and $\mathbf{y}_i(s)= \mathbf{ H}_i(s)\mathbf{u}_i(s) =\mathbf{ L}(s\mathbf{ E} - \mathbf{ A})^{-1} \mathbf{ B}_i\mathbf{u}_i(s)$.

This property can be employed for the parallel computation of the reduced-order model. Each partitioned model of (\ref{super_state}) can be reduced by a projection matrix $\mathbf{X}_i \in \R^{N \times 2k}$ whose columns span the k-dimensional EKS:
\begin{equation*}
\label{Eq:eksm_sup}
\mathcal{K}_k^E(\mathbf{A}_{E},\mathbf{B}_{iE}) = \mathcal{K}_k(\mathbf{A}_{E},\mathbf{B}_{iE}) + \mathcal{K}_k(\mathbf{A}_{E}^{-1},\mathbf{B}_{iE}) =
\end{equation*}
\begin{equation}\label{195}
span \{\mathbf{b}_{iE},  \mathbf{A}_{E}^{-1}\mathbf{B}_{iE}, \mathbf{A}_{E}\mathbf{B}_{iE},\mathbf{A}_{E}^{-2}\mathbf{B}_{iE}, \mathbf{A}_{E}^{2}\mathbf{B}_{iE},\dots,  
\end{equation}
\begin{equation*}
\mathbf{A}_{E}^{-(k-1)}\mathbf{B}_{iE}, \mathbf{A}_{E}^{k-1}\mathbf{B}_{iE}\}
\end{equation*} 
with $\mathbf{B}_{iE} \equiv \mathbf{A}^{-1}\mathbf{B}_i $, and similarly the reduced-order model is obtained by:
\begin{equation} \label{tran_super}
\mathbf{\tilde E}_i = \mathbf{X}_i^T\mathbf{ E}\mathbf{X}_i, \quad \mathbf{\tilde A}_i = \mathbf{X}_i^T\mathbf{ A}\mathbf{X}_i, \quad \mathbf{\tilde B}_i = \mathbf{X}_i^T\mathbf{ B}_i, \quad \mathbf{\tilde L}_i = \mathbf{ L}\mathbf{X}_i
\end{equation}
Moreover, each reduced-order transfer function is computed as:
\begin{equation} \label{h(s)_super}
\mathbf{\tilde H}_i(s)= \mathbf{\tilde L}_i(s\mathbf{\tilde E}_i - \mathbf{\tilde A}_i)^{-1} \mathbf{\tilde B}_i + \mathbf{D}
\end{equation}
Finally, the approximate transfer function of the reduced-order model is computed as:
\begin{equation}\label{H_i}
\mathbf{\tilde H}(s) =   [\mathbf{\tilde H}_1(s),\mathbf{\tilde H}_2(s),\dots,\mathbf{\tilde H}_p(s)]  
\end{equation}

It must be noted that there is no guarantee that the passivity of the reduced-order models obtained using the superposition property is preserved. In recent years, however, the focus of MOR has been shifted from provably passive models to passivity enforcement \textit{after} efficient reduction. A wealth of passivity enforcement techniques, such as \cite{pass}, have been developed to assure passivity of the reduced-order models obtained using the superposition property.

\small
\begin{table*}[!hbt]
	\centering
	\caption{Reduction results of EKS-MM vs MM for industrial IBM power grid benchmarks}
	\label{tab}
	\begin{tabular}{|c|c|c|c|c|c|c|c|c|c|c|}
		\hline
    		\multirow{3}{*}{Ckt} & \multirow{3}{*}{Dimension} & \multirow{3}{*}{\#ports} &
    		
    		\multirow{3}{*}{ROM Order} &
    		\multicolumn{3}{c|}{Moment-Matching (MM)}                            & \multicolumn{4}{c|}{\textbf{EKS Moment-Matching (EKS-MM)} }                  \\ \cline{5-11} 
			&                          &                          &  & \multirow{2}{*}{\#moments}  & \multirow{2}{*}{Max Error} & \multirow{2}{*}{Runtime(s)}      &\multirow{2}{*}{\#moments}& \multirow{2}{*}{Max Error}  &Error Reduction      & \multirow{2}{*}{Runtime(s)} \\ 
		                     &                          &                          &                            &                            &                            &  &            &  &     percentage &                     \\ \hline
		ibmpg1&      44946                    &        600                  & 1200                            &            2               &     0.037                      &   0.146    &    1      &     0.014         &    62.16\%      &                0.146 \\  \hline
		ibmpg2&    127568                      &         500                 &                     2000       &             4               &      0.233                       &  1.206    &     2     &        0.131        &          43.78\%               &  1.277 \\   \hline
		ibmpg3&    852539                      &      800                    &              1600              &              2              &                        0.253     & 11.029     &     1     &        0.146        &          42.29\%               & 11.060 \\  \hline
		ibmpg4&       954545                   &       600                   &                      2400      &          4                  &                        0.233     &   16.642   &      2    &         0.038       &          83.69\%               & 17.981 \\ \hline
		ibmpg5&     1618397                     &      600                    &                  1200          &           2                 &                      0.242       &   10.228   &    1      &          0.063      &          73.97\%               & 10.998 \\ \hline
		ibmpg6&       2506733                   &                  1000        &                   6000         &                  6          &                       0.161      &  19.155    &     3     &         0.130       &         19.25\%               &  21.780 \\ \hline
		ibmpg1t&     54265                     &        400                  &              800              &            2                &                  4.767           &   0.259   &   1       &        1.814        &          61.95\%               & 0.273 \\ \hline
		ibmpg2t&       164897	                   &        800                  &                 3600           &            4                &                       0.785      &  0.250    &      2    &       0.411         &          47.64\%              &  0.268 \\ \hline
	\end{tabular}
\end{table*}

\normalsize
\begin{figure}[!hbt]
    \centering
    \includegraphics[width=1\linewidth,keepaspectratio]{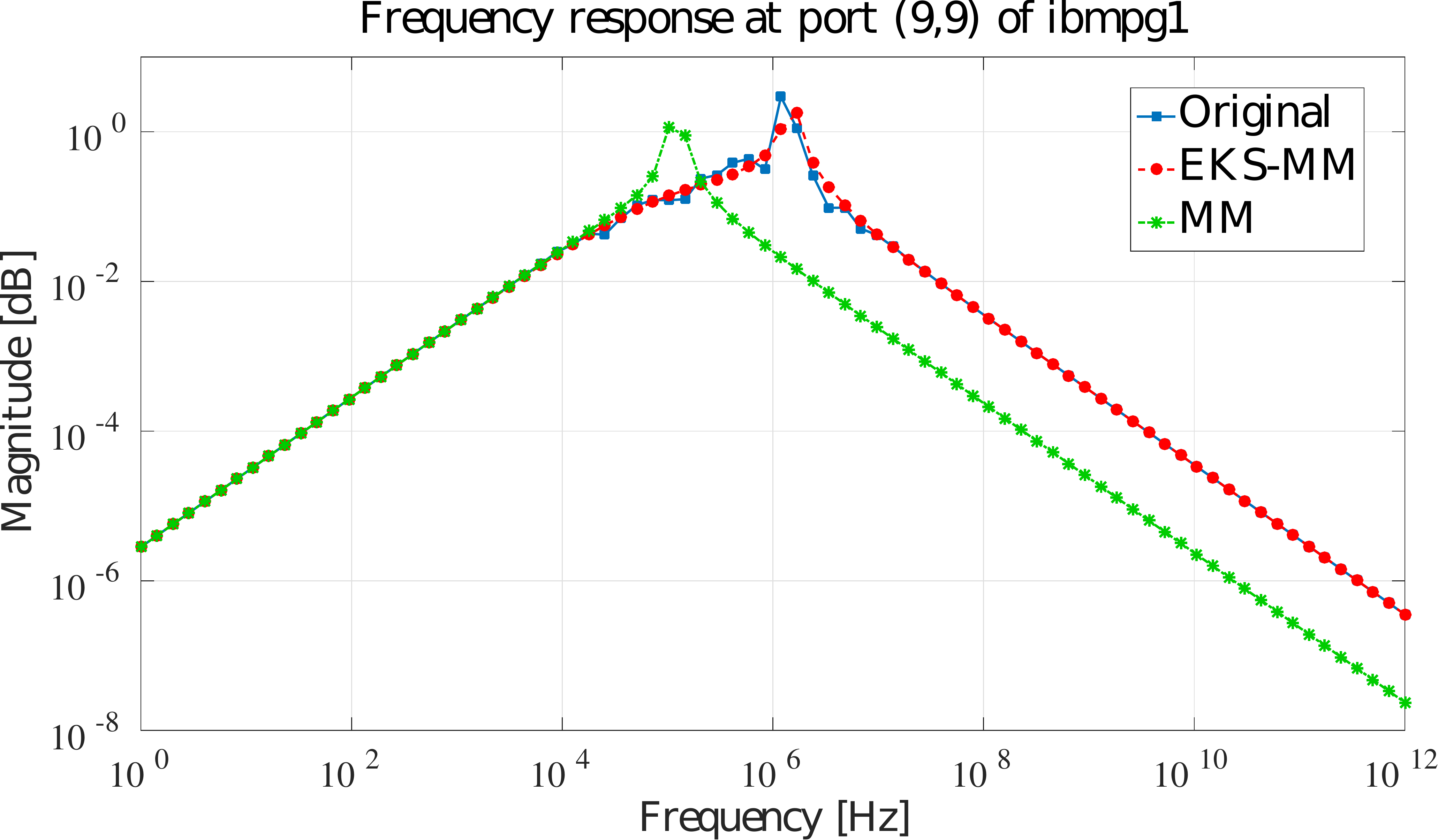}
    
    \vspace{3mm}
    \includegraphics[width=1\linewidth,keepaspectratio]{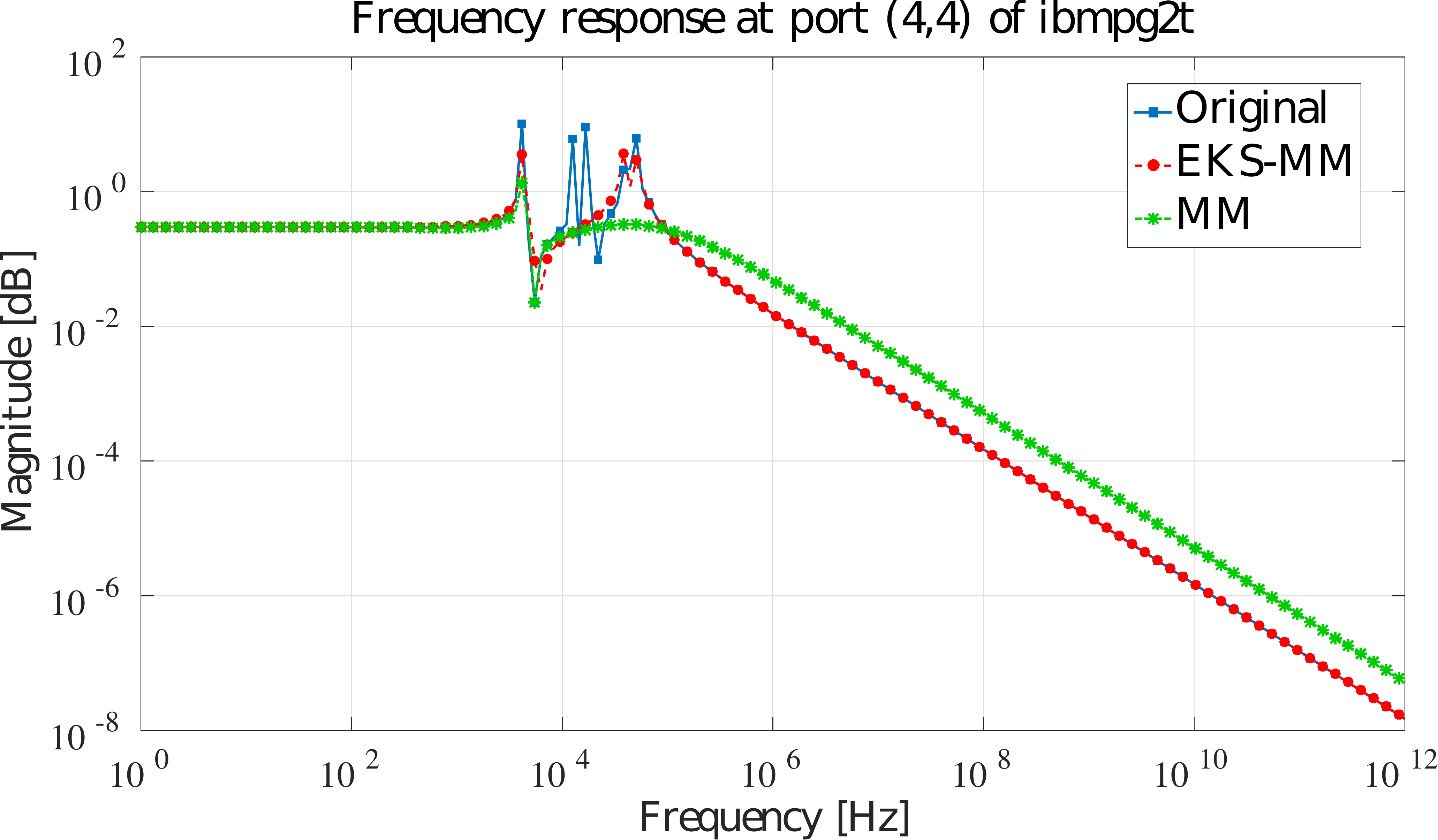}
    \caption{Comparison of transfer functions of ROMs obtained by EKS-MM and MM in the range $[10^0,10^{12}]$, for ibmpg1 and ibmpg2t benchmarks at ports (9,9) and (4,4), respectively.}
    \label{fig:bd1}
\end{figure}



\begin{figure}[!hbt]
    \centering
    \includegraphics[width=1\linewidth,keepaspectratio]{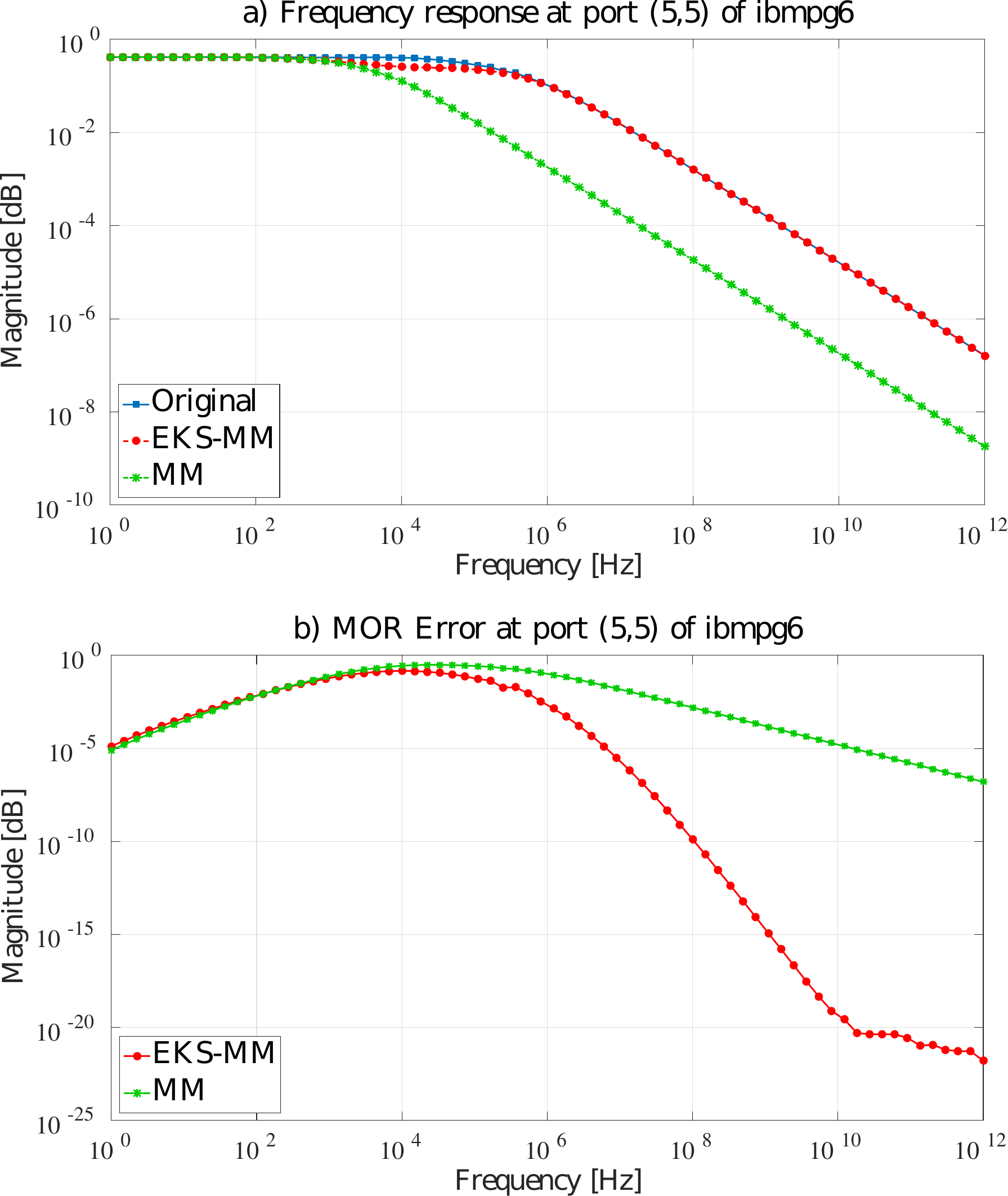}
    \caption{Comparison of transfer functions and absolute error magnitudes of ROMs obtained by EKS-MM and MM in the range $[10^0,10^{12}]$, for ibmpg6 at port (5,5).}
    \label{fig:bd}
\end{figure}
\section{Experimental Results}\label{EXP}
For the experimental evaluation of the proposed methodology, we used the available IBM power grid benchmarks \cite{ibmpg}. Their characteristics are shown in the first three columns of Table \ref{tab}. Note that for the transient analysis benchmarks, ibmpg1t and ibmpg2t, a matrix of energy storage elements (capacitances and inductances) is provided. However, in order to perform transient analysis for the DC analysis bechmarks, ibmpg1 to ibmpg6, we had to add a (typical for power grids) diagonal capacitance matrix with random values on the order of picofarad. In order to evaluate our methodology on singular benchmarks, we enforced the capacitance matrix of ibmpg2 and ibmpg4 to have at least one node that was missing a capacitance connection.
These benchmarks along with ibmpg1t and ibmpg2t were represented as singular descriptor models of (\ref{12}), thus we applied the techniques described in Section \ref{5.3}  for their efficient sparse handling.

EKS-MM was implemented with the procedures described in Section \ref{EKS}, and was compared with a standard MM method also implemented with the superposition property. The reduced-order models (ROMs) were evaluated in the frequency range $[10^0,10^{12}]$ with respect to their accuracy for given ROM order. For our experiments, an appropriate number of matching moments was selected such that the ROM order for both EKS-MM and MM is the same. All experiments were executed on a Linux workstation with a 3.6GHz Intel Core i7 CPU and 32GB memory using MATLAB R2015a. 


The results are reported in the remaining columns of Table \ref{tab},  where \textit{\#moments} refers to the number of moments that matched in order to produce the ROMs, \textit{Max Error} refers to the error between the infinity norms of the transfer functions, i.e. $||\mathbf{\tilde{H}}(s) - \mathbf{H}(s)||_\infty $, \textit{Runtime} refers to the computational time (in seconds) needed to generate each submatrix $\mathbf{H}_i(s)$ of (\ref{H_i}), while \textit{Error Reduction percentage} refers to the error reduction percentage achieved by EKS-MM over MM. It can be clearly verified that, compared to MM for similar ROM order, EKS-MM produces ROMs with significantly smaller error. As depicted in Table~\ref{tab}, the \textit{Error Reduction percentage} ranges from 19.25\% to 83.69\%. The execution time of  EKS-MM is negligibly larger than standard MM for each moment computation, due to the expansion in two points, however the efficient implementation can effectively mask this overhead to a substantial extent and make the procedure applicable to very large circuit models.


To demonstrate the accuracy of our method, we compare the transfer functions of the original model and the ROMs generated by EKS-MM and MM. 
The corresponding transfer functions for one regular (ibmpg1) and one singular (ibmpg2t) benchmark, in the band $[10^0,10^{12}]$, are shown in Fig. \ref{fig:bd1}. Fig. \ref{fig:bd} presents the transfer functions of ROMs produced by EKS-MM and MM along with the absolute errors induced over the original model for a selected benchmark in the same band. As can be seen, the response of EKS-MM ROM is performing very close to the original model,
while the response of MM ROM exhibits a clear deviation. In particular, responses of ROMs produced by MM do not capture effectively the dips and overshoots that arise in some frequencies.
\section{Conclusions}\label{CONCL}
In this paper, we proposed the use of EKS to enhance the accuracy of MM methods for descriptor circuit models. Our method provides clear improvements in reduced-order model accuracy  compared to a standard Krylov subspace MM technique. 
For the implementation,
we made efficient computational choices, as well as adaptations and modifications for large-scale singular models. As a result, the proposed method still remains computationally efficient, introducing only a small overhead in the reduction process.

\end{document}